\documentclass[12pt,a4paper]{article}
\usepackage{amsmath,amsthm,amssymb,amsfonts,hyperref,graphicx,epic,enumerate,float}

\newcommand{\Z}{{\mathbb Z}}
\newcommand{\Nat}{{\mathbb N}}
\newcommand{\gf}{F}
\newcommand{\gfl}{\overline{F}}
\newcommand{\gfc}{\alpha}
\newcommand{\gfcl}{\overline{\alpha}}
\newcommand{\zf}{Z}
\newcommand{\zfl}{\overline{Z}}
\newcommand{\zfc}{A}
\newcommand{\zfcl}{\overline{A}}
\newcommand{\zbl}{\overline{\zb}}
\newcommand{\zb}{\beta} 

\newcommand{\PlanePartitionSet}{\mathfrak{A}}
\newcommand{\StrictPlanePartitionSet}{\mathfrak{B}}
\newcommand{\PBox}{\dashbox{1}(1,1)}
\newcommand{\BOX}[1]{\multiput(0,0)(1,0){#1}{\PBox}}
\newcommand{\Young}{\mathcal{Y}}

\renewcommand{\Re}{\operatorname{Re}}

\theoremstyle{plain}
\newtheorem{thm}{Theorem}
\newtheorem{cor}{Corollary}

\newtheorem{lem}{Lemma}
\newtheorem{prop}{Proposition}

\theoremstyle{definition}
\newtheorem{rem}{Remark}

\newtheorem{defn}{Definition}

\newtheorem{prob}{Problem}

\begin{document}

\author{Johan Andersson\footnote{Department of Mathematics, Stockholm
    University, SE-106 91 Stockholm,
    Sweden. \texttt{johana@math.su.se}}  \and  Jan
  Snellman\footnote{Department of Mathematics, Link\"oping University,   
SE-58183 Link\"oping, Sweden. \texttt{Jan.Snellman@mai.liu.se}}}

\title{On the number of plane partitions and non isomorphic subgroup
  towers of abelian groups} 

\maketitle

\begin{abstract}
  In this paper we study $\gfc_{k,r}(n)$, defined as the number of $k
  \times r$ matrices 
such that $ m_{i,j} \geq  m_{i+1,j} \geq 0$,  $\, \, m_{i,j} \geq
m_{i,j+1}$,  and $m_{1,1}+ \dots +m_{1,r}=n$. We consider the
generating function 
 \begin{gather*}
    \gf_{k,r}(x)=\sum_{n=0}^\infty \gfc_{k,r}(n) x^n. \\ \intertext{We
      use Erhart reciprocity to prove that} 
 \gf_{k,r}(x) = (-1)^{kr} x^{-r(r-1+2k)/2} \gf_{k,r} \left(\frac 1 x \right).
 \end{gather*}
For the special case $k=1$ this result also follows from the classical
theory of partitions,  and for $k=2$ it was proved in
Andersson-Bhowmik \cite{andbhow} with another method. We give an
explicit formula for $\gf_{k,r}(x)$ in terms of Young tableaux. We
then study the corresponding zeta-function \begin{gather}
  \zf_{k,r} (s)=\prod_{p  \text{ prime}} \gf_{k,r}(p^{-s})
  \end{gather}
  and give an application on the average orders of towers of
  abelian groups. In particular we prove that the number of
  isomorphism classes of   
  ``subgroups of subgroups
  of ... ($k-1$ times) ... of abelian groups'' of order at most $N$ is
  asymptotic to $c_k N (\log N)^{k-1}$. This generalises results from
  Erd{\H o}s-Szekeres  \cite{ErdosSzekeres} and  Andersson-Bhowmik
  \cite{andbhow} where the corresponding 
  result was proved for $k=1$ and $k=2$.
  \end{abstract}

\section{Introduction}
From the classical theory of partitions (see e.g. Hardy-Wright
\cite{Hardy} page 281), we have that
if \begin{gather*}
 \gfc_{1,r}(n) = \# \{ q \in \Z^r, 0 \leq q_1 \leq q_i \leq q_{i+1}, \,\, q_1+\dots + q_r=n \}
\end{gather*}
denotes the number of partitions of $n$ into at most $r$ parts, then
\begin{gather*}
  \gf_{1,r}(x)= \sum_{n=0}^\infty \gfc_{1,r}(n) x^n =
   \frac 1 {(1-x) (1-x^2) \cdots (1-x^r)}.
\end{gather*}
In particular this implies that  we have the functional equation
 \begin{gather} \label{funeq1}
   \gf_{1,r}(x) = (-1)^r x^{-r(r+1)/2}  \gf_{1,r} \left( \frac 1 x \right).
 \end{gather}
In the paper  Andersson-Bhowmik \cite{andbhow} the authors considered
the related problem of counting pairs 
$\{p_i,q_i\}_{i=1}^r$ such that  $0 \leq p_i \leq q_i$, $p_i \leq
p_{i+1}$, $q_i \leq q_{i+1}$
  and
 $q_1+\dots+q_r=n$.  A recursion formula was obtained to
 calculate its generation function, and with its help it was proved that if
 \begin{gather*}
    \gf_{2,r} (x)=\sum_{n=0}^\infty \gfc_{2,r} (n) x^n,
 \end{gather*}
 denotes its generating function, then it is a rational
function that satisfies the functional equation
\begin{gather} \label{funeq2}
 \gf_{2,r}(x) = x^{-r(r+3)/2}  \gf_{2,r} \left( \frac 1 x \right).
\end{gather}
We see  that these problems can be viewed in matrix terms as counting
the number of $1 \times r$ integer matrices
\begin{align*}
 &\begin{pmatrix} q_r & \dots & q_1 \end{pmatrix}, \qquad & \qquad
 \qquad &\begin{pmatrix} 0 \leq q_1 \leq \cdots \leq q_r \\ q_1+
   \cdots + q_r=n
\end{pmatrix}  \\
 \intertext{and  $2 \times r$ integer matrices }
  &\begin{pmatrix} q_r & \dots & q_1 \\ p_r & \dots & p_1 \\
  \end{pmatrix} . \qquad & \qquad 
&\begin{pmatrix} 0 \leq p_1 \leq \cdots \leq p_r \\  0 \leq q_1 \leq
\cdots  \leq q_r
\\ p_i \leq q_i, \, \, q_1+ \dots+q_r=n \end{pmatrix}.
\end{align*}
\noindent In this paper we will generalize these problems from
$1\times r$ and $2 \times r$ matrices to $k \times r$ matrices.
\begin{defn} \label{def1}
Let $\gfc_{k,r}(n)$  count the number of $k \times r$ matrices such
that $m_{i,j} \geq  m_{i+1,j}$,  $\, \, m_{i,j} \geq
m_{i,j+1}$, and $m_{1,1}+ \dots +m_{1,r}=n$. Let $\gf_{k,r}(x)$
denote the generating function
\begin{gather}\label{eq:Fdef}
    \gf_{k,r}(x)=\sum_{n=0}^\infty \gfc_{k,r}(n) x^n.
 \end{gather}
\end{defn}
\noindent We prove a functional equation for $\gf_{k,r}(x)$ that
generalizes equations \eqref{funeq1} and \eqref{funeq2}.
\begin{thm}\label{thm:Eproc}
Let $\gf_{k,r}(x)$ be defined by Definition \ref{def1}. Then
\begin{gather*}
  \gf_{k,r}(x) = (-1)^{kr} x^{-r(r-1+2k)/2} \gf_{k,r} \left(\frac 1 x
\right).
\end{gather*}
\end{thm}
 Non-negative integer valued $k \times r$ matrices with
decreasing rows and decreasing columns (as the matrices counted in
Definition \ref{def1}) are also called plane partitions with $k$
columns and $r$ rows. Plane partitions were first studied in
MacMohan \cite{MacMahon} (See also Stanley \cite{Stanley2}). As an
example of a result from the theory: If we define $q_{k,r}(n)$ as
the number of plane partitions with $k$ columns and $r$ rows such
that the sum over all elements in the matrix  equals  $n$, then  the
generating function can be written as
\begin{gather*}
 \sum_{n=1}^\infty q_{k,r}(n) x^n =
 \prod_{i=1}^r \prod_{j=1}^k (1-x^{i+j-1})^{-1}.
\end{gather*}

In our case we define $\gfc_{k,r}(n)$ as the number of plane
partitions with $k$ columns  and $r$ rows such that the sum over the
elements in the first row equals $n$.
In this case the problem will
be more difficult. In fact already for $k=2$ as shown in
Andersson-Bhowmik \cite{andbhow} there seems to be no simple
expression for the generating function.  We also study the limit
case 
\begin{gather} \label{uyew}
    \gfcl_{k}(n)= \lim_{r \to \infty}  \gfc_{k,r}(n), \qquad \text{and}  \qquad   \gfl_{k}(x)= \lim_{r \to \infty}  \gf_{k,r}(x).
\end{gather}
The associated  zeta functions
\begin{gather*}
  \zf_{k,r}(s)= \prod_{p  \text{ prime}} \gf_{k,r} \left(p^{-s}
  \right), \\ \intertext{and}
\zfl_{k}(s)= \prod_{p \text{ prime}}  \gfl_{k}
\left(p^{-s} \right),
\end{gather*}
have interpretations in the context of counting subgroup towers of
abelian groups (of rank at most $r$ or arbitrary rank).
\begin{defn}
  A subgroup tower of a group $G$ of length $k$  is defined
  as a $k$-tuple of groups $(G_1,\ldots,G_k)$ where $G_1=G$ and $G_{j+1}
\subseteq G_{j}$.
\end{defn}
We say that two subgroup towers $G$ and $\tilde G$ are isomorphic if
$G_i \cong \tilde G_i$ for $i=1,\ldots,k$.
 We will use analytic properties of the zeta functions
$\zf_{k,r}(s)$ to prove the following theorem.
\begin{thm} \label{thm2} One has that
\begin{enumerate} \item
 The number of isomorphism classes  of subgroup towers of length
$k$ of abelian  groups of order at most $N$ and rank at most $r$ is
asymptotic to $c_{k,r} N (\log N)^{k-1}$, where \(c_{k,r}\) is a
constant. 
 \item
 The number of  isomorphism classes of subgroup towers of length
$k$ of abelian  groups of order at most $N$ is asymptotic to 
$c_k N (\log N)^{k-1}$.
\end{enumerate}
\end{thm}
This is a classical result of Erd{\H o}s-Szekeres \cite{ErdosSzekeres}
for $k=1$. For $k=2$ it  was proved in 
Andersson-Bhowmik \cite{andbhow}.

\section{The generating function of plane partitions}

\subsection{The functional equation, P-partitions and Ehrhart reciprocity}

\begin{proof}[Proof of Theorem~\ref{thm:Eproc}]
Let \(\PlanePartitionSet_{k,r}\) denote the set of $k \times r$ integer matrices
such that $ m_{i,j} \geq  m_{i+1,j} \geq 0$,  $\, \, m_{i,j} \geq
m_{i,j+1}$, and let \(\StrictPlanePartitionSet_{k,r}\) denote the set obtained by replacing
all inequalities with strict ones. In particular, if 
\((m_{i,j})_{i,j} \in \StrictPlanePartitionSet_{k,r}\) then \(m_{i,j} > 0\).

We give an element \(M=(m_{i,j}) \in \PlanePartitionSet_{k,r}\) weight \(w(M)=\prod
x_{i,j}^{m_{i,j}}\), and introduce the generating functions
\begin{equation}
  \label{eq:functionaleq}
  \begin{split}
    \gf_{k,r}(t_{1,1},\dots,t_{k,r}) &= \sum_{M \in \PlanePartitionSet_{k,r}} w(M) \\
    G_{k,r}(t_{1,1},\dots,t_{k,r}) &= \sum_{M \in \StrictPlanePartitionSet_{k,r}} w(M)
  \end{split}
\end{equation}
Specializing
\begin{equation}\label{eq:spec}
  t_{i,j}=
  \begin{cases}
    x & i = 1 \\
    1 & i > 1
  \end{cases}
\end{equation}
in \(\gf_{k,r}(t_{1,1},\dots,t_{k,r})\) we recover the counting function \eqref{eq:Fdef}.

Denote by \(C_m\) the \(m\)
element chain, and by \(P=P_{k,r}= (C_k \times C_r)\) the \(k\) times
\(r\) ``grid poset''. We have that elements in
\(\PlanePartitionSet_{k,r}\), which are plane partitions, correspond
to \(P\)-partitions, i.e., order-reversing maps from \(P\) to
\(\Nat\), and that elements in \(\StrictPlanePartitionSet_{k,r}\),
which are a special type of plane partitions, correspond
to strict \(P\)-partitions,  i.e., strictly order-reversing maps from \(P\) to
\(\Nat\).
This correspondence is illustrated in Figure~\ref{fig:Ppart}.

\begin{figure}[H]
  \centering
  \setlength{\unitlength}{0.4cm}
  \small{
  \begin{picture}(20,9)
\thicklines
    \multiput(3,0)(4,2){3}{\circle*{0.3}}
    \multiput(0,3)(4,2){3}{\circle*{0.3}}
    \put(3,0){\line(2,1){8}}
    \put(0,3){\line(2,1){8}}
    \multiput(3,0)(4,2){3}{\line(-1,1){3}}

    \put(3.3,-0.6){5}
    \put(7.7,1.5 ){3}
    \put(11.4,3.4){1}
    \put(0.7,2.6){2}
    \put(4.7,4.6){2}
    \put(8.7,6.6){0}

    \multiput(15,5)(2,0){4}{\line(0,1){4}}
    \multiput(15,5)(0,2){3}{\line(1,0){6}}
    \put(15.8,7.7){5}
    \put(17.8,7.7){3}
    \put(19.8,7.7){1}
    \put(15.8,5.7){2}
    \put(17.8,5.7){2}
    \put(19.8,5.7){0}

  \end{picture}
  }
  \caption{P-partitions of the poset \(C_3 \times C_2\) correspond to
    \(2 \times 3\) plane partitions}
  \label{fig:Ppart}
\end{figure}
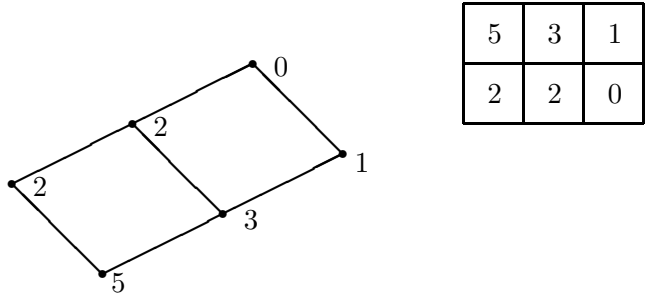

Hence, from the reciprocity theorem for \(P\)-partitions \cite[Thm
4.5.7]{Stanley} (a special case of Ehrhart reciprocity (see also
\cite{Beck}, \cite{Ehrhart})) we have that
\begin{equation}
  \label{eq:reci}
G_{k,r}(t_{1,1},\dots,t_{k,r})  \prod _{i=1}^k\prod_{j=1}^ r t_{i,j}  =
  (-1)^{kr} \gf_{k,r} \left(\frac{1}{t_{1,1}},\dots,\frac{1}{t_{k,r}} \right).
\end{equation}
Furthermore, the obvious bijection
\begin{equation}
  \label{eq:bij}
  \begin{split}
    \phi: \PlanePartitionSet_{k,r} & \to \StrictPlanePartitionSet_{k,r} \\
    M & \mapsto M +
    \begin{pmatrix}
      k+r-1 & k+r-2 & \cdots & k+1 & k \\
      k+r-2 & k+r-3 & \cdots & k & k-1 \\
      \vdots & \vdots & \ddots & \vdots & \vdots \\
      r+1 & r & \cdots & 3 & 2 \\
      r & r-1 & \cdots & 2 & 1
    \end{pmatrix}
  \end{split}
\end{equation}
shows that
\begin{equation}
  \label{eq:Geq}
  G(t_{1,1},\dots,t_{k,r}) =   F  \left(t_{1,1},\dots,t_{k,r} \right)
  \prod_{i=1}^k\prod_{j=1}^r t_{i,j}^{i+j-1}.
\end{equation}
Combining \eqref{eq:reci} and \eqref{eq:Geq} we get that
\begin{equation}
  \label{eq:rel}
  \gf_{k,r}(t_{1,1},\dots,t_{k,r}) = (-1)^{kr}
  \gf_{k,r}\left(\frac{1}{t_{1,1}},\dots,\frac{1}{t_{k,r}} \right)
  \prod_{i=1}^r\prod_{j=1}^k
  t_{i,j}^{-i-j}
\end{equation}
which, using the specialization \eqref{eq:spec} and the relation
\[\sum_{j=1}^r (k+j-1) = r(r-1+2k)/2,\]
 becomes
\begin{equation}
  \label{eq:relsp2}
  \gf_{k,r}(x) = (-1)^{kr}
  \gf_{k,r} \left(\frac{1}{x} \right) x^{-r(r-1+2k)/2}.
\end{equation}
\end{proof}

\subsection{An explicit formula}

Since the total extensions of the poset \(P^*\) is enumerated by
standard Young tableaux of shape \(r^k\), and since a descent
corresponds to a box labeled \(\ell+1\) occurring in a higher row than
the box labeled \(\ell\), we get, by using Theorem 4.5.4 in Stanley \cite{Stanley},
an explicit (but not very efficient)
formula for \(\gf_{k,r}(x)\).

% Namely,
\begin{thm} \label{thm3}
 Let \(T\) be a standard tableaux with shape \(r^k=(r,r,\dots,r)\), let
\(T_\ell\) be the subtableau consisting of the boxes with labels
\(\le \ell\), and let \(c(T_\ell)\) be the number of
boxes in the first row of \(T_\ell\). Let \(d(T_\ell)=c(T_\ell)\) if
\(\ell < rk\) and if in \(T\) the box labeled \(\ell+1\) occurs in a higher row
than the box labeled \(\ell\), and let  \(d(T_\ell)=1\) otherwise. Then
\begin{equation}
  \label{eq:SYT}
  \gf_{k,r}(x) = \sum_{T \in \mathbf{SYT}(r^k)} \prod_{\ell=1}^{kr}
  \frac{t^{d(T_\ell)}}{1-t^{c(T_\ell)}}
\end{equation}
\end{thm}
As an example, we take two rows and three columns. The hook-length
formula \cite{Fulton:Young} shows that there are \(6!/(4*3*3*2*2*1) = 5\)
standard Young tableaux of the desired shape. They are tabulated
below, together with their contribution to \(F_{2,3}(x)\).

\setlength{\unitlength}{0.4cm}

\begin{center} 
  \begin{tabular}{ccc}
     \begin{picture}(3,3)
    \put(0,0){\BOX{3}} \put(0,1){\BOX{3}}
    \put(0.2,1.2){1}     \put(1.2,1.2){2}     \put(2.2,1.2){3}
    \put(0.2,0.2){4}     \put(1.2,0.2){5}     \put(2.2,0.2){6}
  \end{picture} &

     \begin{picture}(3,3)
    \put(0,0){\BOX{3}} \put(0,1){\BOX{3}}
    \put(0.2,1.2){1}     \put(1.2,1.2){2}     \put(2.2,1.2){5}
    \put(0.2,0.2){3}     \put(1.2,0.2){4}     \put(2.2,0.2){6}
  \end{picture} &

     \begin{picture}(3,3)
    \put(0,0){\BOX{3}} \put(0,1){\BOX{3}}
    \put(0.2,1.2){1}     \put(1.2,1.2){2}     \put(2.2,1.2){4}
    \put(0.2,0.2){3}     \put(1.2,0.2){5}     \put(2.2,0.2){6}
     \end{picture}

    \\

    \(\frac{1}{(1-x)(1-x^2)(1-x^3)^4}\) &
    \(\frac{x^2}{(1-x)(1-x^2)^3(1-x^3)^2}\) &
    \(\frac{x^2}{(1-x)(1-x^2)^2(1-x^3)^3}\)

\\

     \begin{picture}(3,3)
    \put(0,0){\BOX{3}} \put(0,1){\BOX{3}}
    \put(0.2,1.2){1}     \put(1.2,1.2){3}     \put(2.2,1.2){5}
    \put(0.2,0.2){2}     \put(1.2,0.2){4}     \put(2.2,0.2){6}
  \end{picture} &

     \begin{picture}(3,3)
    \put(0,0){\BOX{3}} \put(0,1){\BOX{3}}
    \put(0.2,1.2){1}     \put(1.2,1.2){3}     \put(2.2,1.2){4}
    \put(0.2,0.2){2}     \put(1.2,0.2){5}     \put(2.2,0.2){6}
  \end{picture}  &

    \\

    \(\frac{x^3}{(1-x)^2(1-x^2)^2(1-x^3)^2}\) &
    \(\frac{x}{(1-x)^2(1-x^2)(1-x^3)^3}\) &

 \end{tabular}
\end{center}

\noindent For instance, the tableau
\begin{gather*}
 T=  \begin{picture}(3,3)
    \put(0,0){\BOX{3}} \put(0,1){\BOX{3}}
    \put(0.2,1.2){1}     \put(1.2,1.2){2}     \put(2.2,1.2){4}
    \put(0.2,0.2){3}     \put(1.2,0.2){5}     \put(2.2,0.2){6}
  \end{picture} \\ \intertext{contributes}
  \frac{x^2}{(1-x)(1-x^2)^2(1-x^3)^3},
\end{gather*}
as we see by studying the initial subtableux:
\begin{center}
\begin{tabular}{cccccc}
  \begin{picture}(3,3)
   \put(0,1){\BOX{1}}
    \put(0.2,1.2){x}
  \end{picture} &
  \begin{picture}(3,3)
   \put(0,1){\BOX{2}}
    \put(0.2,1.2){x}      \put(1.2,1.2){x}
  \end{picture} &
  \begin{picture}(3,3)
   \put(0,1){\BOX{2}}
   \put(0,0){\BOX{1}}
    \put(0.2,1.2){x}      \put(1.2,1.2){x}  \put (0.2,0.2){1}
  \end{picture} &
  \begin{picture}(3,3)
   \put(0,1){\BOX{3}}
   \put(0,0){\BOX{1}}
    \put(0.2,1.2){x}      \put(1.2,1.2){x}  \put(2.2,1.2){x}
    \put(0.2,0.2){1}
  \end{picture} &
  \begin{picture}(3,3)
   \put(0,1){\BOX{3}}
   \put(0,0){\BOX{2}}
    \put(0.2,1.2){x}      \put(1.2,1.2){x}  \put(2.2,1.2){x}
    \put(0.2,0.2){1}      \put(1.2,0.2){1}
  \end{picture} &
  \begin{picture}(3,3)
   \put(0,1){\BOX{3}}
   \put(0,0){\BOX{3}}
    \put(0.2,1.2){x}      \put(1.2,1.2){x}  \put(2.2,1.2){x}
    \put(0.2,0.2){1}      \put(1.2,0.2){1}     \put(2.2,0.2){1}
  \end{picture}
  \\
   \(\frac{1}{1-x}\) & \(\frac{1}{1-x^2}\) & \(\frac{x^2}{1-x^2}\) &
  \(\frac{1}{1-x^3}\) &   \(\frac{1}{1-x^3}\) &  \(\frac{1}{1-x^3}\)
\end{tabular}
\end{center}
By adding the five terms corresponding to the standard tableaux we get
\begin{gather}
 \gf_{2,3}(x)=\frac{\left( 1 + 2 x + 2 x^2 + x^3 + x^4 \right)  
     \left( 1 + x + 2 x^2 + 2 x^3 + x^4 \right) }{{\left( -1 + x \right) }^
      6 {\left( 1 + x \right) }^3 {\left( 1 + x + x^2 \right) }^4},
\end{gather}
which coincide with the result calculated in Andersson-Bhowmik \cite{andbhow}. By similar reasoning we obtain
\begin{gather}
 \gf_{3,2}(x)= \frac{x^4+2x^3+4x^2+2x+1}{(x-1)^6 (x+1)^5}, \\ \intertext{and}
 \gf_{4,2}(x)=\frac{x^6+3x^5+9x^4+9x^3+9x^2+3x+1}{(x-1)^8 (x+1)^7}. 
\end{gather}
We see that Theorem \ref{thm3} implies that the generating function is a rational function.
\begin{cor} \label{cor1a}
 The generating function $\gf_{k,r}(x)$ can be written as
 \begin{gather*}
    \gf_{k,r}(x)= \frac{Q_{k,r}(x)} { \left((1-x) \cdots (1-x^{kr}) \right)^{kr} }
 \end{gather*}
 where $Q_{k,r}(x)$ is a polynomial of degree $(k-1) r (kr+r-1)$
\end{cor}
\begin{proof} This follows from Theorem \ref{thm3}. The degree of the polynomial $Q_{k,r}(x)$ follows from Theorem \ref{thm:Eproc}. \end{proof}

\begin{rem}
 The estimation for the degree of the denominator in Corollary \ref{cor1a} is an overestimation. In fact for e.g. $k=2$ it can be shown (either by the recursion formula or by more careful consideration of the Young tableaux)  that $\gf_{2,r}(x)=P_{2,r}(x)/((1-x)^{2} (1-x)^{3} \cdots (1-x^r)^{r+1})$ where $P_{2,r}(x)$ is a polynomial of degree $r(-5+3r+2 r^2)/6$.
\end{rem}

\subsection{Further properties of the generating function}
We calculate the first few values of $\gfc_{k,r}(n):$
\begin{lem} \label{uu}
 One has that \begin{enumerate}[(i)]
      \item  $\gfc_{1,r}(n)=p_r(n)$.
        \item  $\gfc_{k,r}(0)=1$.
       \item  $\gfc_{k,r}(1)=k$.
       \item $\displaystyle \gfc_{k,r}(2)= \begin{cases} k(k+1)/2, &
           r=1, \\ k(k+1), & r \geq 2, \end{cases} $
       \item $\displaystyle \gfc_{k,1}(n)= \binom{n+k-1}{k-1}$.
       \end{enumerate}
   \end{lem}
\begin{proof}
\begin{enumerate}[(i)]
  \item This is just the classical restricted partition function (Hardy-Wright \cite{Hardy} page 281).
  \item The only matrix that will contribute is the zero matrix.
  \item The matrices that will contribute  have the first $j$ rows
$(1,0,\ldots,0)$ and
        $k-j$ rows $(0,\ldots,0)$, for $j=0,\ldots,k-1$. There are $k$ such matrices.
  \item The matrices that will contribute will either have
        \begin{enumerate}
          \item $j_1$ rows $(2,0,\ldots,0)$, $j_2$ rows $(1,0,\ldots,0,1)$ and $j_3$ rows $(0,\ldots,0)$ such that $j_1+j_2+j_3=k$ and  $j_1 \geq 1$, or
      \item $j_1$ rows $(1,1,0,\ldots,0)$, $j_2$ rows $(1,0,\ldots,0)$ and $j_3$ rows $(0,\ldots,0)$ such that $j_1+j_2+j_3=k$ and  $j_1 \geq 1$.
         \end{enumerate}
         The number in each case will  be $\binom {k+1} 2$ and in
         general we get the contribution
        $2 \, \binom {k+1} 2=k(k+1)$. If $r=1$ only the first case
        will contribute and we will instead get just $\binom {k+1} 2$.
  \item The matrices that will contribute will have the rows $(a_i)$ for $i=1,\ldots,k$ and
    $n=a_1 \geq a_2 \geq \cdots \geq a_{k} \geq 0$. There are $\binom{n+k-1}{k-1}$ such matrices.
\end{enumerate}
\end{proof}
We see that Lemma~\ref{uu} $(v)$ implies that
\begin{gather} \label{ttr}
 \gf_{k,1}(x)= (1-x)^{-k}.
\end{gather}

One may ask the following question: What happens with the generating
function $\gf_{k,r}(x)$
when $r$ or $k$ tends to infinity. From Lemma~\ref{uu} $(iii)$ it
follows that  $\lim_{k
  \to \infty} \gf_{k,r}(x)$ is divergent. However in the case when $r
\to \infty$ we can take the limit.
We define
\begin{gather} \label{gfxkdef}
   \gfcl_k(m)=\gfc_{k,m}(m).
\end{gather}

\begin{lem} \label{uuuu} One has that
\begin{gather} \notag
  0 \leq \gfc_{k,r}(m) \leq \gfc_{k,r+1}(m), \\ \intertext{and in particular}
    \label{uty}  \gfc_{k,r}(m) = \gfcl_k(m). \qquad \qquad (r \geq m)·
\end{gather}
\end{lem}
 \begin{proof}
     It is clear that $\gfc_{k,r}(m) \geq 0$, since it is a counting
     function. That it increases in  $r$ follows from the fact that
     every $k\times r$ matrix which is counted in  $\gfc_{k,r}(m)$ will
     correspond to the $k \times (r+1)$ matrix where we adjoin a zero
     column as the last column, which is counted in
     $\gfc_{k,r+1}(m)$. That $\gfc_{k,r}(m) = \gfcl_k(m)$ for $(r
     \geq m)$ follows from the fact that
      we can have at most $m$ non
     zero columns under the given conditions (the maximum number of
     non zero columns will be attained exactly when the first row has
     $m$ ones and $(r-m)$ zeroes.
 \end{proof}
    \begin{lem} \label{uuuu2} One has that
\begin{gather} \notag
   q(m) \leq \gfcl_k(m) \leq ((m+1)q(m))^k,
\end{gather}
 where $q(m)$ denote the classical partition function.
\end{lem}
\begin{proof}
     The lower bound is obtained by counting the matrices with the
     first row an arbitrary
     partition of $n$ and $k-1$ rows identically zero.
    For the upper bound we use the fact that
     each row in a matrix that we count for $\gfc_{k,m}(m)$ will be a
     classical partition for some number $0 \leq j \leq m$. Hence we
     have the inequality
     \begin{gather*}
          \gfc_{k,m}(m) \leq  \left(q(0)+\cdots+q(m) \right)^k
     \end{gather*}
     Since the classical partition function is an increasing function
     this implies that
     \begin{gather*}
          \gfc_{k,m}(m) \leq   ((m+1)q(m))^k.
     \end{gather*}
\end{proof}
We introduce the generating function
\begin{gather}
  \gfl_k(x)=\sum_{n=1}^\infty  \gfcl_k (n) x^n.
\end{gather}
and we  prove (The proof is essentially  the same as the proof of
Lemma 2 of Andersson-Bhowmik \cite{andbhow}):
\begin{lem} \label{u13}
 With $\gf_{k,r}(x)$ and $\gfl_k(x)$ defined as above one has that $\gf_{k,r}(x)$ and $\gfl_k(x)$ are analytic
 functions in the unit disc with integer power series coefficients
 such that $\gf_{k,r}(0)=\gfl_k(0)=1$.
Furthermore the function $\gfl_k(x)$
satisfies the inequality
 \begin{gather*}
    \gfl_k(x) \geq \frac 1 {\prod_{k=1}^\infty (1-x^k)}. \qquad \qquad (0<x<1)
 \end{gather*}
\end{lem}

\begin{proof}
The power series coefficients of $\gf_{k,r}(x)$ and $\gfl_k(x)$ are integers since they are counting functions and
by Lemma  \ref{uu}  $(ii)$ and eq. \eqref{uty}  we have that
$\gfc_{k,r}(0)=1$ and $\gfl_k(0)=1$,
 which implies $\gf_{k,r}(0)=\gfl_k(0)=1$.
By the well known generating function for the classical partition function
 \begin{gather} \label{o112}
  \sum_{n=0}^\infty q(n) x^n= \frac 1 {\prod_{n=1}^\infty (1-x^n)}, \qquad \qquad (0<x<1)
  \end{gather}
 and the lower bound in Lemma \ref{uuuu2}
 \begin{gather*}
   q(n) \leq \gfcl_k(n),
 \end{gather*}
this gives us the lower bound in Lemma \ref{u13}. Equation \eqref{o112} also
implies that the
generating function of the partition function is analytic in the unit disc,
and hence the
classical partition function $q(n)$ is of subexponential order. This implies
that   $((n+1)q(n))^k$ is
of subexponential order and by the upper bound in Lemma \ref{uuuu2}, so is
$\gfc_{k,r}(n)$, and also $\gfcl_k(n)$ since $0 \leq \gfc_{k,r}(n)\leq \gfcl_k(n)$.
This  proves that $\gfl_k(x)$ and $\gf_{k,r}(x)$ are analytic in the unit disc.
\end{proof}

\subsection{The polynomials \(\gfc_{k,r}(n)\)}

\begin{prop} \label{prop1}
  For fixed \(n,r\), the quantity \(\gfc_{k,r}(n)\) is a polynomial of degree $n$ in
  \(k\), with leading coefficient \(\gfc_{1,r}(n)/n!\).

  Similarly, for fixed \(n\), the quantity \(\gfcl_{k}(n)\) is a polynomial of degree $n$ in
  \(k\), with leading coefficient \(\gfcl_{1}(n)/n!\).
\end{prop}
\begin{proof}
  A plane partition \(M=(m_{ij})\) of dimension \(k \times r\) can be regarded as a
  sequence of \(k\) partitions \(m_1 \le m_2 \le \dots \le m_k\) with
  at most \(r\) parts,
  where each \(m_i\) corresponds to the \(k+1-i\)'th row of \(M\). Hence,
  \(M\) can be viewed as a   multichain (i.e., a chain with possible
  repetitions) of length \(k\) in the 
  restricted Young lattice \(\Young_r\) of partitions with at most  \(r\)
  parts. If the last row of \(M\) sums to \(n\), then the associated
  multichain ends at level \(n\) in the ranked lattice \(\Young_r\).

  Now consider the principal order ideal \(I\) in \(\Young_r\) generated by
  \(m_k\). The whole multichain \(m_1 \le m_2 \le \dots \le m_k\) is
  contained in \(I\), and the number of such \(k\)-multichains in
  \(I\) is given by \(Z(I;\, k) = \zeta^k(\emptyset,m_k)\), where
  \(\zeta\) is the (combinatorial) zeta function of \(\Young_r\). It is
  well known that this expression is a polynomial in \(k\); indeed, it
  is called the Zeta-polynomial of \(I\), see for instance
  \cite[IV:2]{Aigner:Combinatorial}.

  Clearly, we get the desired quantity \(\gfc_{k,r}(m)\) by summing over all
  Zeta-polynomials of principal ideals of partitions of \(n\) with at
  most \(r\) parts,
  \begin{equation}
    \label{eq:idealSumRest}
    \gfc_{k,r}(n) = \sum_{\substack{\lambda \vdash n\\ \ell(\lambda)
        \le r}} Z(I_{\le \lambda}; \, k)
  \end{equation}
  This is a finite sum of polynomials in \(k\), hence a polynomial in
  \(k\).

It follows from standard properties of the Zeta-polynomial (see again
\cite[IV:2]{Aigner:Combinatorial})  that each
\(Z(I_{\le \lambda})\) is an integer-valued polynomial in \(k\) of degree \(n\) with
non-negative coefficients, hence that \(\gfc_{k,r}(n)\) and
\(\gfcl_{k}(n) = \gfc_{k,n}(n)\) both have degree \(n\). Evaluated at 1, this
polynomial is \(\gfc_{1,r}(n)\), the number of partitions of \(n\)
with at most \(r\) parts. We get from the elementary theory of
integer-valued polynomials (see for instance \cite[Corollary
1.3]{Stanley:CombCom}) that the leading coefficient of the polynomial
\(\gfc_{k,r}(n)\) is \(\gfc_{1,r}(n)/n!\).

To obtain the result for the unrestricted coefficient
\(\gfcl_{k}(n)\), replace \(\Young_r\) with \(\Young\).
\end{proof}

 Furthermore:
\begin{prop}
  The polynomial \(\gfc_{k,r}(n)\in \mathbf{Q}[k]\) is divisible by
  \((k+s)\) for all \(s < \frac{n}{r+1}\). 
  There is some \(c\), independent of \(n\), such that the polynomial \(\gfcl_{k}(n)\in
  \mathbf{Q}[k]\) is divisible by
  \((k+s)\) for all integers \(s < c \sqrt{n}\).
\end{prop}
\begin{proof}
  Since \(\Young_r\) is a locally finite distributive lattice, every
  interval is a finite distributive lattice. Hence, the M\"obius function
  for an interval \([\lambda,\tau]\) is either 0 or
  \((-1)^{\textrm{rk}(\lambda) -\textrm{rk}(\tau)}\). Thus, we can
  apply a result of Stanley's \cite[Proposition 4.9]{Aigner:Combinatorial} which says that
  (putting \(P=I_{\le \lambda} \subset \Young_r\) for some partition
  of \(n\) with at
  most \(r\) parts)
  \begin{equation}
    \label{eq:negvals}
    Z(P; \, -k) = (-1)^{\textrm{rk}(P)} \overline{Z}(P; \, k)
  \end{equation}
  where \(\overline{Z}(P; \, k)\) counts the number of \(k\)-multichains
  \begin{equation}
    \label{eq:neqchains}
  0 \le z_1 \le \cdots \le z_k =1
  \quad
  \text{ with }
  \quad
  \mu(z_{i-1},z_i) \neq 0
  \quad
  \text{ for all i.}
  \end{equation}
  Furthermore, \(\mu(z_{i-1},z_i) \neq 0\) if and only if
  \([z_{i-1},z_1]\) is a Boolean lattice, which happens if and only if
  the supremum of the points in the interval is a Boolean algebra.

  If there is an upper bound \(L\) on the length of Boolean
  subintervals of \(P\), then it follows that \(\mu(\lambda,\tau)=0\) whenever
  \(\textrm{rk}(\lambda) -\textrm{rk}(\tau) > L\). Thus, the
  minimal length of a chain \eqref{eq:neqchains} is \(\frac{N}{L}\),
  where \(N\) is the length of \(P\), i.e., the length of the longest
  chain in \(P\). Hence  \(\overline{Z}(P; \, k)=0\) for \(k < \frac{L}{L}\), hence,
  by \eqref{eq:negvals},
  the same holds for \(Z(P; \,k)\). It follows that the polynomial
  \(Z(P; \,k)\) is divisible by \(k+s\) for \(s <  \frac{N}{L}\).

  A partition of \(s\) with \(\le r\)
  parts can be covered by at most \(r+1\) partitions, the supremum of
  which has rank \(\le s+r+1\), so the maximum length of a Boolean
  subinterval in \( I_{\le \lambda} \subset \Young_r\) is \(r+1\). By
  the above, it follows that the polynomial
  \(Z(P; \,k)\) is divisible by \(k+s\) for \(s <  \frac{n}{r+1}\).

  We now turn to partitions with an unlimited number of parts.
Let the truncated Young lattice \(\Young^{\le n}\) consist of
partitions of \(s \le n\).
To estimate, by the above method,  how many zeroes at negative
integers the polynomial
\(\gfcl_k(n)\) will have, we would need to bound the size of the
intervals in  \(\Young^{\le n}\) that are Boolean algebras. Such an
interval would look like \([\lambda,\tau]\) with
\(\lambda=(\lambda_1,\dots,\lambda_v)\), and \(\tau\) the supremum of
the elements covering \(\lambda\). A partition \(\tilde{\lambda}\)
covering \(\lambda\)   will either be
\(\tilde{\lambda}=(\lambda_1,\dots,\lambda_v,1)\) or
else
\(\tilde{\lambda}=(\lambda_1,\dots, \lambda_i+1,\dots,\lambda_v,1)\).
The position \(i\) were the extra element is inserted must either be
1, or else \(\lambda_i < \lambda_{i-1}\). If there are \(\ell\) such
positions (the position at the end included) then there are \(\ell\)
elements covering \(\lambda\), the supremum \(\tau\) of these elements is a
partition of \(\lvert \lambda \rvert + \ell\), so \([\lambda,\tau]\)
is a Boolean algebra of length \(\ell\).

The partition \(\lambda = (s,s-1,\dots,1)\) is a partition of
\(\binom{s+1}{2}\) and is covered by \(s+1\) partitions; the supremum is a
partition of \(\binom{s+1}{2} + s+1= \binom{s+2}{2}\), which should be
no larger than \(n\) for the interval to fit inside \(\Young^{\le n}\).

This is  maximal,
so that any Boolean algebra inside \(\Young^{\le n}\) have length less
or equal to
\(\frac{\sqrt{1+8n} -3}{2}\). The second assertion now follows.
\end{proof}

Using \eqref{eq:idealSumRest} and a MAPLE package by Stembridge
\cite{Stembridge:Posets}, we can calculate \(\gfc_{k,r}(n)\) using the
following simple commands:
\small{
\begin{verbatim}
with(SF);
with(posets);
read("young_lattice");
alphaPart := proc(part,varname)
    zeta(young_lattice(part),varname);
end;
alphaPartlist := proc(PLIST,varname)
    local q;
    add(alphaPart(q,varname), q = PLIST);
end;
aRN := proc(r,n,varname)
    local L,q;
    L := select(q->nops(q) <= r, combinat[partition](n));
    L := map(q->sort(q,`>`),L);
    alphaPartlist(L,varname);
end;

> factor(aRN(4,4,k));
                                                2
                             5 k (k + 2) (k + 1)
                             --------------------
                                      12
\end{verbatim}
}

Recall that \begin{math}
\gfc_{k,1}(n) = \binom{k+n-1}{n}
\end{math}.
We tabulate the polynomials \(\gfcl_{k}(n)\), \(\gfc_{2,k}(n)\) and
\(\gfc_{2,k}(n)\) below.

\begin{align*}
%  \label{eq:gfctab1}
\gfcl_{k}(1) &= k
\\
\gfcl_{k}(2) &= k \left( k+1 \right)
\\
\gfcl_{k}(3) &= 1/6\,k  \left( k+1 \right) \left( 4\,k+5 \right)
\\
\gfcl_{k}(4) &= {\frac {5}{12}}\,k \left( k+1 \right)^{2} \left( k+2 \right)
\\
\gfcl_{k}(5) &= {\frac {1}{60}}\,k \left( k+1 \right) \left( k+2 \right)    \left( 13\,
{k}^{2}+36\,k+21 \right)
\\
\gfcl_{k}(6) &= {\frac {1}{180}}\,k \left( k+1 \right)  \left( k+2
\right)^{2} \left( 19\,{k}^{2}+58\,k+33
 \right)
\\
\gfcl_{k}(7) &= {\frac {1}{2520}}\,k \left( k+
1 \right) \left( k+2 \right) \left( k+3 \right)      \left( 116\,{k}^{3}+508\,{k}^{2}+688\,k+263 \right)
\\
\gfcl_{k}(8) &= {\frac {1}{10080}}\,k \left( k
+1 \right) \left( k+2 \right) \left( k+3 \right)      \left(
191\,{k}^{4}+1338\,{k}^{3}+3297\,{k}^{2}+3330\,k+
1084 \right)
\intertext{and}
    \gfc_{k,2}(1) &= k
\\
\gfc_{k,2}(2) &= k \left( k+1 \right)
\\
\gfc_{k,2}(3) &= 1/2\,k \left( k+1 \right) ^{2}
\\
\gfc_{k,2}(4) &= 1/4\,k \left( k+2 \right)  \left( k+1 \right) ^{2}
\\
\gfc_{k,2}(5) &= 1/12\,k \left( k+2 \right) ^{2} \left( k+1 \right) ^{2}
\\
\gfc_{k,2}(6) &= 1/36\,k \left( k+3 \right)  \left( k+2 \right) ^{2} \left( k+1
 \right) ^{2}
\intertext{and}
    \gfc_{k,3}(1) &= k
\\
\gfc_{k,3}(2) &= k \left( k+1 \right)
\\
\gfc_{k,3}(3) &= 1/6\,k \left( 4\,k+5 \right)  \left( k+1 \right)
\\
\gfc_{k,3}(4) &= 1/24\,k \left( 9\,k+7 \right)  \left( k+2 \right)  \left( k+1 \right)
\\
\gfc_{k,3}(5) &= {\frac {1}{120}}\,k \left( k+2 \right)  \left( k+1 \right)  \left( 21
\,{k}^{2}+52\,k+27 \right)
\\
\gfc_{k,3}(6) &= {\frac {1}{240}}\,k \left( k+2 \right)  \left( k+1 \right)  \left( 17
\,{k}^{3}+82\,{k}^{2}+125\,k+56 \right)
\end{align*}

\subsection{The growth of the coefficients $\gfcl_k(n)$}
For $k=1$ we have that $\gfcl_1(n)$ equals the classical partition
function $q(n)$. Thus $\gfcl_{k}(n)$ is a proper generalisation of the
partition function $q(n)$. For $q(n)$ we have good asymptotics by a
theorem of Hardy-Ramanujan \cite{HardyRamanujan}
\begin{gather} \label{rr}
 q(n) \sim \frac{e^{\pi \sqrt{2n/3}}}{4n\sqrt{3}}.
\end{gather}
This is a strong result and it seems difficult to obtain a
similar formula for the general case. Equation \eqref{rr} implies that
\begin{gather} \label{rr23}
   \frac{\log q(n)}{\sqrt n}  = \pi \sqrt{2/3}+ o(1).
\end{gather}
We will thus study
\begin{gather} \label{rr24}
   \frac{\log \gfcl_{k}(n)}{\sqrt n}.
\end{gather}
We improve on the lower
bound in Lemma \ref{uuuu2}.
\begin{prop}
 One has that $\gfcl_{k}(kn) \geq q(n)^k$.
\end{prop}
\begin{proof}
Let $r=kn$. And let $B=\{q_{i,j}\}$ be a $k \times r$ matrix such that
$q_{i,j} \geq q_{i,j+1} \geq 0$ and the sum of each row
$q_{i,1}+\cdots+q_{i,r} = n$ is a partition of $n$.
 It is clear that there are exactly $q(n)^k$ such matrices. For
each matrix $B$ of this type we can construct a matrix $A$
\begin{gather*}
A= \begin{pmatrix} q_{1,j}+\cdots+q_{1,k}& \cdots &
    q_{r,1}+\cdots+q_{r,k}  & 0 & \cdots & 0 \\
   q_{1,j}+\cdots+a_{1,k-1}& \cdots &
    q_{r,1}+\cdots+q_{r,k-1}  & 0 & \cdots & 0 \\
       \vdots &   & \vdots & \vdots & &  \vdots \\
    q_{1,j}+q_{1,2}& \cdots &
    q_{r,1}+q_{r,2}  & 0 & \cdots & 0 \\
  q_{1,j}& \cdots &
    q_{r,1} & 0 & \cdots & 0
   \end{pmatrix}
\end{gather*}
This matrix will be counted in $ \alpha_{k,kn}(kn)$. Hence
$\gfcl_{k}(kn)= \gfc_{k,kn}(kn) \geq q(n)^k$.
 \end{proof}

Together with the upper bound in Lemma  \ref{uuuu2} and eq. \eqref{rr23}
this implies the following Corollary:
\begin{cor} \label{cor99} Suppose that $k \geq 1$ is an integer. Then
 \begin{gather*}
     \sqrt k + o(1) \leq \frac{\log \gfcl_k(n)}{\pi \sqrt {2n/3}}
     \leq  k + o(1). 
\end{gather*}
\end{cor}

\begin{figure}
% \hfill
\hskip -24 pt
 \begin{minipage}[t]{8.2 cm}
\includegraphics[width=8 cm, height=4 cm]{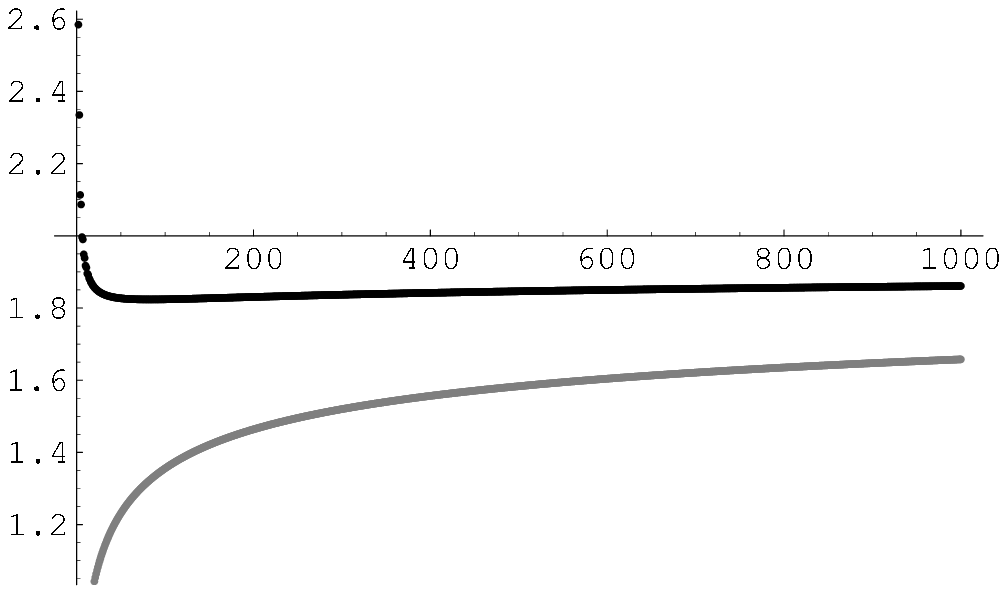}
\vskip -10 pt
\end{minipage}
\begin{minipage}[t]{5cm}
\includegraphics[width=5.5 cm, height=4 cm]{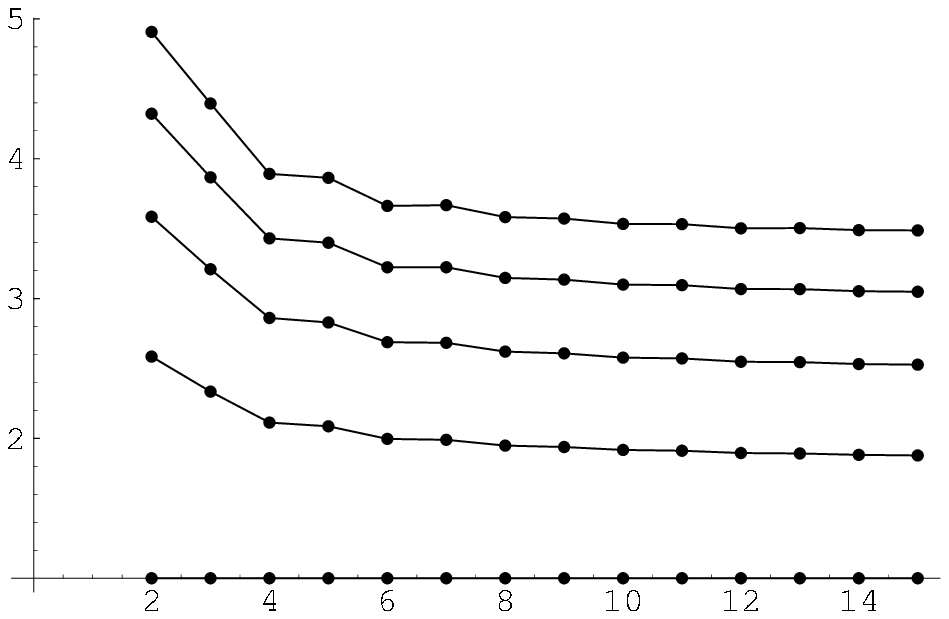}
\end{minipage}
\caption{Left graph: black
  line  $\frac{\log \gfcl_{2}(n)}{\log q(n)}$,  grey line  $\frac{\log \gfcl_{2}(n)}{\pi
    \sqrt{2n/3}}$.
  Right graph: from top to bottom, $\frac{\log \gfcl_{k}(n)}{\log q(n)}$
   for $k=5,\ldots,1$.
}
\label{figa} \end{figure}
The case $k=2$ can be studied numerically by means of the recursion
formula given as Proposition 2 in
Andersson-Bhowmik \cite{andbhow}. For $k \geq 3$ we can use the $\gfc_{k,r}(n)$-polynomials that we
calculated,  although the algorithm is less efficient than
the recursion formula for the case $k=2$ and it  will be difficult to calculate
$\gfc_{k,r}(n)$ for $n \geq 30$.  We show related plots in
Figure~\ref{figa} above.

Even though Corollary \ref{cor99}  is somewhat of an improvement to the trivial lower bound
$1+o(1)$ we would like to know better. The graphs in Figure~\ref{figa} suggests that equation \eqref{rr24} might have a limit and we propose the following problem.
 \begin{prob}
  Find an asymptotic formula for $\log \gfcl_k(n)$.
\end{prob}

\subsection{Nonisomorphic subgroup towers of abelian $p$-groups}
Let $p$ be a prime. An abelian $p$-group is an abelian group of order
$p^n$. Each abelian $p$-group of rank at most $r$ and order $p^n$ is isomorphic  to
a group
\begin{gather}
  G= \Z/(p^{q_1}\Z)\oplus  \cdots \oplus  \Z/(p^{q_r}\Z) \qquad (0 \leq q_1 \leq \cdots \leq
  q_r,  \, \, q_1+ \cdots + q_r=n)
\end{gather}
With the definition of subgroup towers (Definition \ref{def1}) we see
that each subgroup tower of length $k$ and maximal 
group of order $p^n$ and rank $r$ will be isomorphic to
\begin{gather*}
  (G_1,\ldots,G_k)\cong
  \left(\oplus_{j=1}^r\Z/(p^{m_{1,j}}\Z),\ldots,\oplus_{j=1}^r
    \Z/(p^{m_{k,j}}\Z) \right), \\ \intertext{such that} 
       m_{i,j} \geq m_{i+1,j} \geq 0, \qquad m_{i,j} \geq m_{i,j+1}
      \qquad \text{and} \qquad m_{1,1}+ \cdots+m_{1,r}=n.
\end{gather*}
These are exactly our plane partitions with $k$ rows and $r$ columns
such that the sums of the 
elements in the first row equals $n$ that we already studied. Hence,
we get the following Lemma:
\begin{lem} \label{uy} 
  \begin{enumerate}[(i)]
  \item The number of isomorphism classes of  subgroup towers of
    length $k$ such that 
  the maximal group has order $p^n$ and rank at most $r$ equals
  $\gfc_{k,r}(n)$.
  \item The number of  isomorphism classes of  subgroup towers of
    length $k$ such that 
  the maximal group has order $p^n$ equals
  $\gfcl_k(n)$.
  \end{enumerate}
\end{lem}
We will see how this will give average orders for non isomorphic
subgroup towers in the next section.

\section{The zeta function and nonisomorphic  subgroup towers of
  finite abelian groups}

\subsection{Average orders}

If $G$ is a finite abelian group of order $n$ we have by the
fundamental theorem of finite abelian groups that
\begin{gather*}
 G  \cong  G_{p_1^{a_1}} \oplus \cdots \oplus  G_{p_r^{a_r}}
\end{gather*}
where $n=p_1^{a_1} \cdots p_m^{a_m}$ and $G_{p_j^{a_j}}$ is a $p$-group of order
$p^{a_j}$. By Lemma \ref{uy} we see that
\begin{prop} \label{uy23} Suppose that
    $n=p_1^{a_1} \cdots p_m^{a_m}$. Then
 \begin{enumerate}[(i)]
            \item
  The number of  isomorphism classes of  subgroup towers of length $k$
  such that 
  the maximal group has order $n$ and rank at most $r$ equals
  \[\zfc_{k,r}(n)=\prod_{j=1}^m  \gfc_{k,r}(a_j).\]
  \item
  The number of isomorphism classes of   subgroup towers of length $k$
  such that 
  the maximal group has order $n$ equals
  $\zfcl_{k}(n)=\prod_{j=1}^m  \gfcl_{k}(a_j)$. \end{enumerate}
\end{prop}
We see that  $\zfc_{k,r}(n)$ and $\zfcl_k(n)$ are multiplicative
functions and as in Andersson-Bhowmik  \cite{andbhow} we can introduce the zeta
functions
\begin{gather}
    \zf_{k,r}(s)=\prod_{p \text { prime}} \gf_{k,r}(p^{-s})= \sum_{n=1}^\infty \zfc_{k,r} (n) n^{-s}, \\ \intertext{and}
    \zfl_k (s)=\prod_{p \text { prime}} \gfl_k(p^{-s})= \sum_{n=1}^\infty \zfcl_k(n) n^{-s}.
  \end{gather}
and they will have interpretations in terms of nonisomorphic subgroup
towers of abelian groups. By Lemma \ref{u13} and Dahlquist's theorem
\cite{D} we obtain the Proposition.

\begin{prop} \label{oj} Let $\epsilon>0$. There exist a positive
  integer $P$ such that 
 \begin{gather*}
   \zf_{k,r}(s) = \zeta(s)^k \times \left( \prod_{\substack{p<P \\ p\text{ prime}}}
    \gf_{k,r}(p^{-s}) (1-p^{-s})^k \right) \times  \left(
    \prod_{m=2}^\infty \zeta_P(ms)^{\zb_{k,r}(m)} \right), \\
    \intertext{and}  \zfl_k(s) = \zeta(s)^k \times \left(
      \prod_{\substack{p<P \\ p\text{ prime}}} 
        \gfl_k(p^{-s})(1-p^{-s})^k \right) \times \left(
        \prod_{m=2}^\infty \zeta_P(ms)^{\zb_k(m)} \right), \\ 
    \intertext{valid for $\Re(s)>\epsilon$, where}
        \zeta_P(s) = \zeta(s) \times \left( \prod_{\substack{p<P \\ p\text{ prime}}}
          (1-p^{-s})  \right) \,  = \, \, 
   \prod_{\substack{p  \geq P\\ p  \text{ prime}}} (1-p^{-s})^{-1}. \\ \intertext{Furthermore}
\zb_{r,k}(m)=  \sum_{d|m} \mu \left(\frac m d \right) \frac d m
B_{r,k}(d), \text{ where }
  \log F_{k,r}(x)  = \sum_{m=1}^{\infty} B_{k,r}(m)x^{m},
\end{gather*}
$\zbl_{k}(m)=\zb_{k,m}(m)$, and  $\zb_{k,r}(m)$ and $\zbl_{k}(m)$ are integers.
\end{prop}
\begin{proof} This follows from Lemma  \ref{u13} and  Dahlquist's   \cite{D}
  Lemma 2. \end{proof}

From this, the  fact that $F_r(x)$ has no zeroes for $0<x<1$ (positive
power series coefficients), and from the explicit values of
$\gfc_{k,r}(1)$, $\gfc_{k,r}(2)$ 
given by Lemma \ref{uu} and \eqref{ttr}, the following Corollary follows:
\begin{cor} \label{cor3}
  One has that
   \begin{enumerate}[(i)]
      \item $\displaystyle \zf_{k,r}(s) = 
   \begin{cases} \zeta(s)^k, & r=1, \\ \zeta(s)^k \zeta(2s)^{k(k+1)/2}
     G_{k,r}(s), & r \geq 2, \end{cases}$  
        \item $\displaystyle \zfl_{k}(s) = 
        \zeta(s)^k \zeta(2s)^{k(k+1)/2} \overline{G}_{k}(s),$
   \end{enumerate}
     where $G_{k,r}(s)$ and $\overline{G}_k(s)$ are Dirichlet series
     absolutely convergent and 
    without real zeroes for
    $\Re(s)>1/3$.
\end{cor}
The average order of the Dirichlet series coefficients
$\zfc_{k,r}(n)$ and $\zfcl_{k}(n)$ which count the relevant subgroup
towers (Lemma \ref{uy23}) will come from the pole of the corresponding
zeta-functions at $s=1$ and by a standard Tauberian argument
\cite[Theorem 4.20]{
Grunewald},
Corollary \ref{cor3} implies 
Theorem \ref{thm2}.

\subsection{The polynomials $\zb_{k,r}(n)$ and analytic properties  of the zeta-function}
By  the inequality in Lemma \ref{u13}
 \begin{gather*}
    \gfl_k(x) \geq \frac 1 {\prod_{k=1}^\infty (1-x^k)} \qquad \qquad (0<x<1)
 \end{gather*}
 it is clear that $\gfl_k(x)$ can not be written as a finite product
\begin{gather*}
  \gfl_k(x)=\prod_{j=1}^m (1-x^j)^{b_j}. \qquad \qquad (b_j \in \Z)
\end{gather*}
Hence Dahlquist's theorem also implies the following proposition.
\begin{prop}
  The zeta-functions $\zfl_k(s)$ can be meromorphically continued to
  $\Re(s)>0$ but not beyond the imaginary axis.
\end{prop}
This problem can also be studied for $\zf_{k,r}(s)$.
For further analytic information about the zeta functions $\zfl_k(x)$
and $\zf_{k,r}(s)$ we need to study the coefficients $\zbl_{k}(n)$ and
$\zb_{k,r}(n)$. By their definition in Proposition \ref{oj} and the fact that
$\gfc_{k,r}(n)$ are polynomials in $k$ of degree $n$ (Proposition \ref{prop1}) the following Proposition follows.
\begin{prop}
    For fixed \(n,r\), the quantity \(\zb_{k,r}(n)\) is a polynomial in
  \(k\), as is \(\zbl_{k}(n)\).
\end{prop}
We tabulate the first few polynomials.
\begin{equation}
  \label{eq:gfctab33}
  \begin{split}
\zbl_{k}(1) &= k \\
\zbl_{k}(2) &= \frac{1}{2} k (k+1) \\
\zbl_{k}(3) &= \frac{1}{6} k (k+1) (k+2) \\
\zbl_{k}(4) &=-\frac{1}{12} (k-3) k (k+1) (k+2) \\
\zbl_{k}(5) &= \frac{1}{120} k (k+1) (2 k+1) \left(k^2+k+18\right) \\
\zbl_{k}(6) &= \frac{1}{120} k (k+1) (k+2) \left(k^3-6 k^2-4
  k+29\right) \\
\zbl_{k}(7) &= -\frac{(k-3) k (k+1) (k+2) \left(8 k^3+49
    k+48\right)}{1260}
\end{split} \end{equation}
In Andersson-Bhowmik we calculated the first values for $k=2$
\begin{gather*}
  \zbl_2(1),\ldots,\zbl_2(15)=
  2, 3, 3,4 ,2 ,6 , 1, 4, 6, 2, 0, 12, -1, -2, 9, \ldots \\ \intertext{For $k=3,4$ we
    calculate}
  \zbl_3(1),\ldots,\zbl_3(15) = 3, 6, 6, 10, 0, 21, -5, 0, 51,
  -42, -6, 110, -100, 151, -492 \\
 \zbl_4(1),\ldots,\zbl_4(13) = 4, 10, 10, 20, -10, 57, -19, -72, 324, -370, -92, 1137, -2406
\end{gather*}
 In general we see that the polynomials $\zbl_{k}(1), \zbl_{k}(2)$ and
 $\zbl_{k}(3)$ are positive for $k \geq 1$ and the polynomial $\zbl_{k}(4)$
 is negative for $k \geq 4$. This implies some analytical properties
 of the zeta function.
\begin{prop}
  One has that $\zfl_k(s)$  and $\zf_{k,r}(s)$ for $r \geq 4$
 have no poles for $\Re(s)>1/4$ except for a
  pole of order $k$ at $s=1$, a pole of order $k (k+1)/2$ at $s=1/2$
  and a pole of order $k(k+1)(k+2)/6$ at $s=1/3$. Under the Riemann
  hypothesis it follows  that $\zfl_k(s)$ has no poles for
  $\Re(s)>1/8$ except for possible poles at $1/4,1/5,1/6,1/7$.
\end{prop}

\providecommand{\bysame}{\leavevmode\hbox to3em{\hrulefill}\thinspace}
\def\cprime{$'$} \def\cprime{$'$} \def\cprime{$'$} \def\cprime{$'$}

\bibliographystyle{alphaurl}

\end{document}